\documentclass[english,a4paper,11pt]{article}

\usepackage{amssymb,amsmath,amsfonts,amsthm,epsfig,pstricks,graphics,tikz}

\usepackage[paper=a4paper,left=30mm,right=30mm,top=25mm,bottom =25mm]{geometry}

\usepackage[width=0.85\textwidth ]{caption}

\newtheorem{theorem}{Theorem}[section]

\usepackage{hyperref}

\providecommand{\keywords}[1]
{
  \small	
  \textbf{\textit{Keywords---}} #1
}

\title{\LARGE{\textbf{{Integer Sequences from Circle Divisions in Rational Billiards}}}}

\author{Daniel Jaud}
\date{}
\begin{document}
\maketitle
\flushbottom

\hrulefill\\ %Zeile mit Linie füllen

\begin{abstract}
\noindent
We study rational circular billiards. By viewing the trajectory formed after each reflection point to another inside the circle as the number of circle divisions into regions we derive a general formula for the number of division regions after each reflection. This will give rise to an integer division sequence. Restricting to the special cases $\vartheta =\frac{q}{2q+1}\cdot 2\pi$ we show that the number of regions after each reflection $n$ is beautifully related to Gauss 's arithmetic series.
\end{abstract}

\vspace*{0.4cm}
\keywords{rational billiards, circle division, integer sequence, arithmetic series}

\vspace*{0.7cm}
\hrulefill\\ %Zeile mit Linie füllen
\tableofcontents
%\hrulefill\\
%\newpage
\hrulefill\\

\section{Introduction and setup}\label{sec:Introduction}
\begin{figure}[htb]
\centering
\begin{tikzpicture}[scale=1.5]
%\fill[white] (-3,0) circle (0.05);
\draw (0,0) circle (2);
\fill[black] (0,0) circle (0.05);
\fill[black] (0,-2) circle (0.1) node[below] at (0,-2.2) {$P_0$};
\draw[dashed] (210:2)--(270:2)--(330:2)--(30:2);
\draw(270:2)--(330:2)--(0,0)--cycle;

\draw[dashed] (-2,-2)--(2,-2) node[right] {circle tangent};
\node[left] at (-4,0) {};
\draw (0,-0.75) arc(270:330:0.75) node[above] at (0.25,-0.7) {$\vartheta$};
\draw (0,-1.25) arc (90:30:0.75) node[left] at (0.5,-1.6) {$\alpha$};

\draw (-1,-2) arc (180:150:1) node[right] at (-0.95,-1.8) {$\beta$};

\draw (1,-2) arc (0:30:1) node[left] at (0.95,-1.8) {$\beta$};

\end{tikzpicture}

\caption{Graphical representation of the billiard trajectory, the incident angle $\alpha$ and the angle $\vartheta$ between two consecutive reflection points. \label{fig:setup}}
\end{figure}
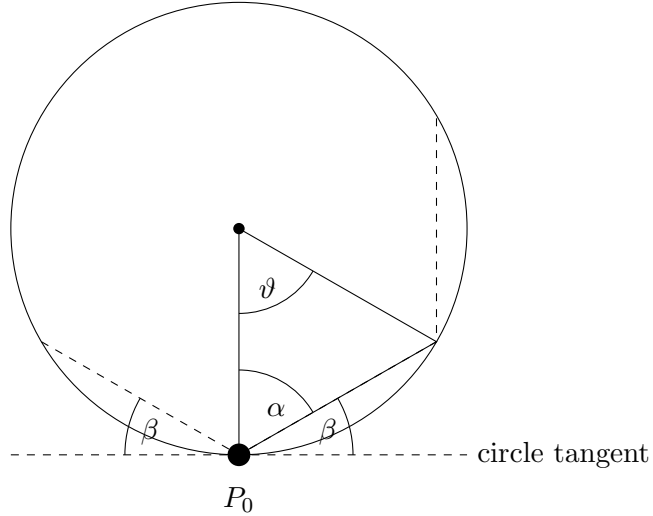

The simplest realization of a billiard table is given by a circular region $\cal C$. Here we are interested in the trajectory of the billiard ball which is considered to be a point particle moving inside the boundary of the circle with unit velocity and bouncing elastically along the boundary $\partial \cal C$. Due to the elastic collision the incident and reflected angle $\alpha$ with respect to the normal to the boundary are identical. Further, because of rotational symmetry the system is fully described by the angle $\vartheta$ made by two consecutive scattering points with the circle \cite{billiard1, billiard2} (see figure \ref{fig:setup}).

From the figure we determine $\vartheta$, depending on the incident angle $\alpha$, to be
\begin{equation}
\vartheta =\pi -2\alpha.
\end{equation}

The natural question of this problem is for which values of $\vartheta$ closed, i.e. periodic, orbits exist. This means how do we need to choose the value of $\vartheta$ such that after a given number of collisions $n$ the billiard ball returns to its initial position $P_0$.
Obviously, due to rotational symmetry this case applies if
\begin{align}
q\cdot \vartheta &= 2\pi \cdot p\\
\leftrightarrow \vartheta &=\frac{p}{q}\cdot 2\pi,
\end{align}

where $q,p\in\mathbb{N}$ such that $gcd(p,q)=1$ and $\frac{p}{q}< \frac{1}{2}$\footnote{The cases $p/q >1/2$ simply correspond to a billiard ball moving in the opposite direction, i.e. clockwise. Since this gives no new insights to the systems behavior we restrict to the cases smaller or equal $1/2$.}. The result tells us that if $\vartheta$ is a rational multiple of $2\pi$ we always end up with periodic orbits \cite{billiard1}. $q$ determines the number of periods that the circle is surrounded before returning to $P_0$. 

\begin{figure}[htb]
\centering
\begin{tikzpicture}[scale=1.5]
\draw (0,0) circle (2);
\foreach \i in {0,1,2,3,4} \draw (360*2/5*\i:2)--({360*2/5*(\i+1)}:2);

\foreach \i in {0}  \fill[black] (360*1/5*\i:2) circle (0.07) node[right] {$P_{\i}$};
\foreach \i in {1}  \fill[black] (360*1/5*\i:2) circle (0.07) node[above] {$P_{\i}$};
\foreach \i in {2,3}  \fill[black] (360*1/5*\i:2) circle (0.07) node[left] {$P_{\i}$};
\foreach \i in {4}  \fill[black] (360*1/5*\i:2) circle (0.07) node[below] {$P_{\i}$};
\draw[blue,dashed] (2,0)--(0,0)-+(72:2);
\draw[blue] (0.5,0) arc (0:72:0.5) node[above] at (0.25,-0.05) {$\vartheta /p$};
\draw[red,->] (2.5,0) arc (0:144:2.5) node[above] at (0,2.7) {$(p-1)$ star corners between };
\end{tikzpicture}

\caption{Example for a $5-$star formed by the trajectory corresponding to $\vartheta =\frac{2}{5}\cdot 2\pi$. \label{fig:setup2}}
\end{figure}
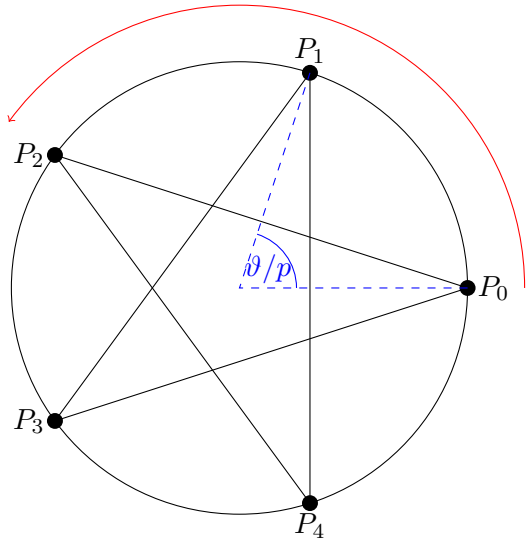

The trajectory forms in general regular star shaped figures consisting of $q$ corners equally distributed along $\partial {\cal C}$\footnote{If $p=1$ the trajectory describes regular polygons. Star shaped figures arise for $p\geq 2$.}. Two consecutive star corners $P_i$ and $P_{i+1}$ ($i\in\{0;1;\dots ;q-1\}$) lying on $\partial \cal C$ form an angle of $\vartheta/p$ with the circle (see figure \ref{fig:setup2}). Thus the parameter $p$ tells us that from one scattering point to the next $p-1$ corner points lie in between. If $n$ is the number of scattering then the billiard trajectory follows the path formed by the star corners 
\begin{equation}
P_i=P_{p\cdot n~\mbox{\small mod}~ q},
\end{equation}
where in polar coordinates (circle radius $R$) the positions of $P_i$ are given by
\begin{equation}
P_i\left(R\cos\left(\frac{\vartheta}{p}\cdot i\cdot 2\pi\right),R\sin\left(\frac{\vartheta}{p}\cdot i\cdot 2\pi\right)\right).
\end{equation}

In the following sections we are interested in deriving a general formula for the number of regions the circle is divided by the billiard trajectory after each scattering $n$. Thereby, we continue as follows: In section 2 we are going to show that any two billiard trajectories intersect only in a single point, i.e. it doesn't occur that three or more paths cross in a single point. This is a rather technical side remark to the main results that we will develop in section 3. In particular, we will derive a general formula for the number of circle divisions after each scattering. The final section 4 will be dealing with the special cases when 
\begin{equation}
\vartheta =\frac{p}{2p+1}\cdot 2\pi~~~~~p\in\mathbb{N}.
\end{equation}
In these cases the general formula of section 3 will be significantly simplified containing the well known Gau\ss ~formula for the integer series
\begin{equation}
\sum_{i=1}^n i =\frac{n\cdot(n+1)}{2}.
\end{equation}

\section{Number of intersecting lines}

\begin{theorem}
In rational circular billiards there only intersect two lines formed by the billiard trajectory in a single point.
\end{theorem}

\begin{proof}

For the proof we consider the rotational angle $\vartheta=\frac{p}{q}\cdot 2\pi$, where $gcd(p,q)=1$ and $\frac{p}{q}<\frac{1}{2}$. The trajectory forms a regular star shaped geometry with reflection points $P_i$ on the circle are determined by the angle $\Theta_i=i\cdot \frac{\vartheta}{p},~i\in \{0;\dots;p-1\}$, i.e. $\Theta_i =\angle P_0MP_i$. We thus can consider the following setup of possible intersecting lines (see figure \ref{fig:intersection}) where without loss of generality due to the rotational symmetry we can restrict to the case when the billiard ball starts at the point $P_0$ characterized by $\Theta_0$.
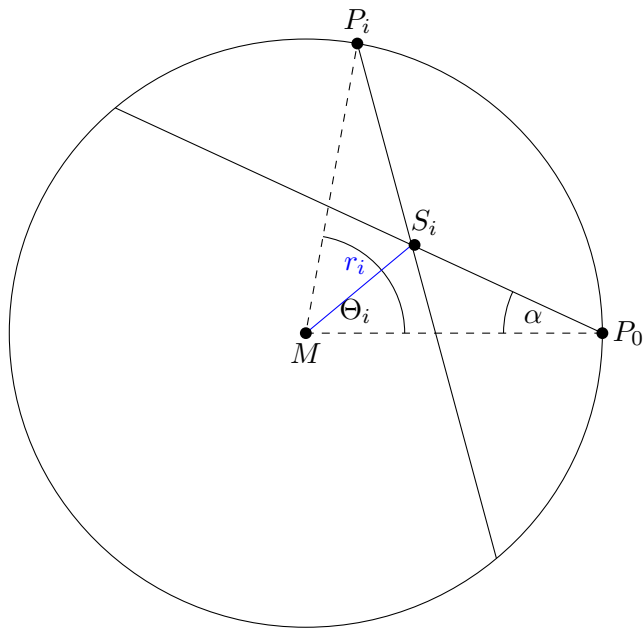
\begin{figure}[htb]
\centering
\begin{tikzpicture}[scale=1.3]
\draw (0,0) circle (3);
\draw[dashed] (3,0)--(0,0)-+(80:3);
\draw (1,0) arc (0:80:1) node[above] at (0.5,0) {$\Theta_i$};
\draw (3,0)-+(130:3);
\draw (80:3)-+(310:3);
\draw[color=blue] (0,0)-+(40:1.4) node[above] at (0.5,0.5) {$r_i$};
\draw (2,0) arc (180:155:1) node[above] at (2.3,0) {$\alpha$}; 
\fill[black] (1.1,0.9) circle (0.06) node[above] at (1.2,0.9) {$S_i$};
\fill[black] (80:3) circle (0.06) node[above] {$P_i$};
\fill[black] (3,0) circle (0.06) node[right] {$P_0$};
\fill[black] (0,0) circle (0.06) node[below] {$M$};

\end{tikzpicture}

\caption{Possible intersecting trajectory lines. \label{fig:intersection}}
\end{figure}
The angle $\alpha$, i.e. the scattering angle with respect to the circle normal, is given by $\displaystyle \alpha=\frac{\pi-\vartheta}{2}$. Making use of the law of sines one obtains for the distance $r_i$ of the intersection point $S_i$ with respect to the circle center $M$ (radius of the circle $R$):
\begin{equation}
r_i=R\cdot \frac{\sin(\alpha)}{\sin(\pi-\alpha - \frac{\Theta_i}{2})}=R\cdot \frac{\cos(\frac{p}{q}\pi)}{\cos(\frac{p}{q}\pi-\frac{i}{q}{\pi} )}.
\end{equation}
Since holds
\begin{equation}
0<\frac{p}{q}\pi-\frac{i}{q}{\pi}<\frac{\pi}{2} ~~~\forall i\in\{0;\dots ;p-1\}
\end{equation}
the $\cos(\frac{p}{q}\pi-\frac{i}{q}{\pi} )$ in the denominator is strictly monotonic increasing. This implies that $r_i$ is strictly monotonic decreasing resulting in the fact that there exist no intersection points $S_i$ resulting from intersecting more than two trajectory lines.
\end{proof}

\begin{figure}[htb]
\begin{minipage}{0.45\textwidth}
\centering
\begin{tikzpicture}[scale=1.5]
\draw (0,0) circle (2);
\foreach \i in {0,1,2,3,4,5,6,7,8,9,10,11,12,13} \draw (360*3/7*\i:2)--({360*3/7*(\i+1)}:2);

%\draw[blue] (0,0) circle (0.7137);
%\draw[red] (0,0) circle (0.4939);
\end{tikzpicture}
\end{minipage}
\hfill
\begin{minipage}{0.45\textwidth}
\centering
\begin{tikzpicture}[scale=1.5]
\draw (0,0) circle (2);
\foreach \i in {0,1,2,3,4,5,6,7,8,9,10,11,12,13} \draw (360*3/7*\i:2)--({360*3/7*(\i+1)}:2);

\draw[blue] (0,0) circle (0.7137);
\draw[red] (0,0) circle (0.4939);
\node[right] at (2.1,0) {${\cal C}_0$};
\node[right,blue] at (0.72,0) {${\cal C}_1$};
\end{tikzpicture}
\end{minipage}
\caption{\textit{left:} Star shaped formed by billiards with $\vartheta=\frac{3}{7}\cdot 2\pi$. \textit{right:} Intersection points of the trajectories aligned on circles of radii $r_i$. \label{fig:circleintersections}}
\end{figure}
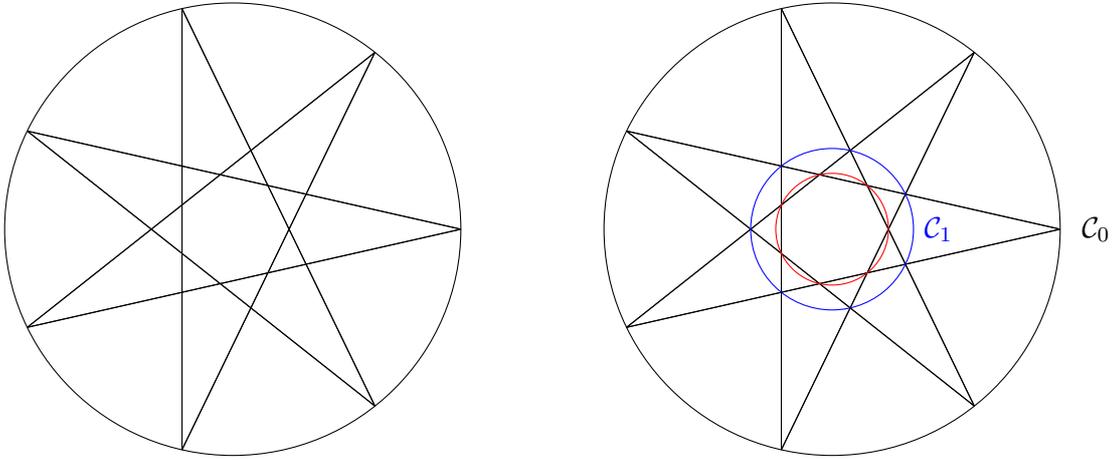

As a direct consequence of the proof it results that the intersection points of the trajectories (intersections along $\partial \cal C$ not counted) again lie on $p-1$ different circles of radii $r_i$ (see figur \ref{fig:circleintersections}).

Interestingly, each subregion ${\cal C}_i$ made by smaller circles of radii $r_i$ again give rise to valid periodic orbits where the angle $\vartheta$ with respect to ${\cal C}_i$ changes according to
\begin{equation}
\vartheta\rightarrow \vartheta_i=\frac{p-i}{q} 2\pi~~~~~~i\in\{0;\dots ;p-1\}.
\end{equation} 

See, for example, the blue circle in figure \ref{fig:circleintersections} corresponding in the example to ${\cal C}_1$.

As a direct consequence of the above shown result for the radii $r_i$ we see that when the ball returns to its initial point $P_0$ the circle has been divided into $p\cdot q+1$ regions.

\begin{proof}
Let $\partial {\cal C}_i$ be the circle of radius $r_i$ with $i\in\{0;\dots ,p-1\}$. For $i=p-1$, i.e. the smallest circle, $\partial {\cal C}_{p-1}$ is the boundary for a regular polygon of $q$ corners, i.e. enclosing one region. Each additional circle $\partial {\cal C}_i$ for $0<i<p-2$ is the boundary of a star shaped figure adding $q$ additional regions (see the yellow regions in figure \ref{fig:regions} for $\partial {\cal C}_{p-2}$) for each circle except for $i=0$ where we have the $q$ regions formed by the star's corners and $q$ regions from the outer circle $\partial {\cal C}_0$ with the next inner circle $\partial {\cal C}_1$ (see green and blue regions in figure \ref{fig:regions}). Thus the total number of regions $f_{total}$ is given by

\begin{equation}
f_{total}=1+\left(\sum_{i=1}^{p-2}q\right)+2q=pq+1.  \label{eq:areas} 
\end{equation}

\begin{figure}[htb]
\centering
\begin{tikzpicture}[scale=1.5]
\fill[yellow] (0,0) circle (1);
\fill[white,draw=white] (0:0.4939)--(0,0)-+(360/7:0.4939);
\fill[white,draw=white] (360/7:0.4939)--(0,0)-+(2*360/7:0.4939);
\fill[white,draw=white] (2*360/7:0.4939)--(0,0)-+(3*360/7:0.4939);
\fill[white,draw=white] (3*360/7:0.4939)--(0,0)-+(4*360/7:0.4939);
\fill[white,draw=white] (4*360/7:0.4939)--(0,0)-+(5*360/7:0.4939);
\fill[white,draw=white] (5*360/7:0.4939)--(0,0)-+(6*360/7:0.4939);
\fill[white,draw=white] (6*360/7:0.4939)--(0,0)-+(7*360/7:0.4939);

\draw (0,0) circle (2);
\fill[blue!70!white,opacity=1] ({(360/7)*0}:2)-+({(360/7)*0 + (360/14)}:0.7137)-+({(360/7)*0 +360/7}:2) arc ({(360/7)*0 +360/7}:{(360/7)*0}:2);
\fill[blue!70!white,opacity=1] ({(360/7)*1}:2)-+({(360/7)*1 + (360/14)}:0.7137)-+({(360/7)*1 +360/7}:2) arc ({(360/7)*1 +360/7}:{(360/7)*1}:2);
\fill[blue!70!white,opacity=1] ({(360/7)*2}:2)-+({(360/7)*2 + (360/14)}:0.7137)-+({(360/7)*2 +360/7}:2) arc ({(360/7)*2 +360/7}:{(360/7)*2}:2);
\fill[blue!70!white,opacity=1] ({(360/7)*3}:2)-+({(360/7)*3 + (360/14)}:0.7137)-+({(360/7)*3 +360/7}:2) arc ({(360/7)*3 +360/7}:{(360/7)*3}:2);
\fill[blue!70!white,opacity=1] ({(360/7)*4}:2)-+({(360/7)*4 + (360/14)}:0.7137)-+({(360/7)*4 +360/7}:2) arc ({(360/7)*4 +360/7}:{(360/7)*4}:2);
\fill[blue!70!white,opacity=1] ({(360/7)*5}:2)-+({(360/7)*5 + (360/14)}:0.7137)-+({(360/7)*5 +360/7}:2) arc ({(360/7)*5 +360/7}:{(360/7)*5}:2);
\fill[blue!70!white,opacity=1] ({(360/7)*6}:2)-+({(360/7)*6 + (360/14)}:0.7137)-+({(360/7)*6 +360/7}:2) arc ({(360/7)*6 +360/7}:{(360/7)*6}:2);

\fill[green!70!white,opacity=1] (360/7*0:0.4939)-+(360/7*0 +360/14:0.7137)-+(360/7*0:2)-+(360/7*0-360/14:0.7137)--cycle;
\fill[green!70!white,opacity=1] (360/7*1:0.4939)-+(360/7*1 +360/14:0.7137)-+(360/7*1:2)-+(360/7*1-360/14:0.7137)--cycle;
\fill[green!70!white,opacity=1] (360/7*2:0.4939)-+(360/7*2 +360/14:0.7137)-+(360/7*2:2)-+(360/7*2-360/14:0.7137)--cycle;
\fill[green!70!white,opacity=1] (360/7*3:0.4939)-+(360/7*3 +360/14:0.7137)-+(360/7*3:2)-+(360/7*3-360/14:0.7137)--cycle;
\fill[green!70!white,opacity=1] (360/7*4:0.4939)-+(360/7*4 +360/14:0.7137)-+(360/7*4:2)-+(360/7*4-360/14:0.7137)--cycle;
\fill[green!70!white,opacity=1] (360/7*5:0.4939)-+(360/7*5 +360/14:0.7137)-+(360/7*5:2)-+(360/7*5-360/14:0.7137)--cycle;
\fill[green!70!white,opacity=1] (360/7*6:0.4939)-+(360/7*6 +360/14:0.7137)-+(360/7*6:2)-+(360/7*6-360/14:0.7137)--cycle;

\foreach \i in {0,1,2,3,4,5,6,7,8,9,10,11,12,13} \draw (360*3/7*\i:2)--({360*3/7*(\i+1)}:2);

\draw[blue,ultra thick] (0,0) circle (0.7137);
\draw[red, ultra thick] (0,0) circle (0.4939);
 
\end{tikzpicture}

\caption{Graphical representation of enclosed regions by each circle.} \label{fig:regions}
\end{figure}
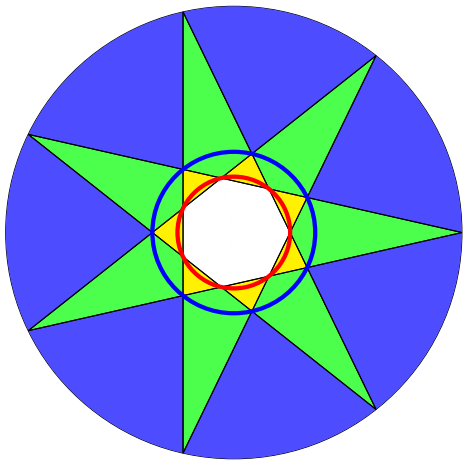

\end{proof}

Equation \eqref{eq:areas} is in agreement with Euler's formula \cite{Euler} for planar graphs
\begin{equation}
f=1+e-v
\end{equation}
where we excluded the region outside of the circle (thus number 1 instead of 2), $f$ is the number of faces, $e$ the number of edges and $v$ the number of vertices, i.e. in our case all line intersections including the reflection points in the circle. In our setup, i.e. when the billiard ball has performed a whole periodic orbit, it holds
\begin{align}
e&=2pq,\\
v&=pq,
\end{align}
and thus
\begin{equation}
f=f_{total}=2pq-pq+1=pq+1.
\end{equation}

\section{General formula for circle division sequence}
In this section we derive a general formula for the division of the circle in regions by the rational billiard trajectory with $\vartheta =\frac{p}{q}\cdot 2\pi$. It turns out that it is useful to rewrite $q$ in terms of $p$ as
\begin{equation}
q=mp+r,
\end{equation}
with $m=\lfloor \frac{q}{p}\rfloor \in\mathbb{N} $ and the rest term $r\in\{1;2;\dots ; p-1\}$. With this definition we can now turn to the circle division starting with the circular region ${\cal C}$ and a billiard ball starting from the initial positions $P_0$. The initial number of divided regions such is $f_0=1$. For each consecutive scattering $i\geq 1$ the circle is divided in one additional region as long as no whole circle surrounding is performed, i.e. the maximum number for $i_{max}$ reflections is given by
\begin{align}
i_{max,1}\cdot \vartheta &<2\pi \\
\leftrightarrow i_{max,1}&< \frac{2\pi}{\vartheta}=\frac{q}{p}=\frac{mp+r}{p}=m+\frac{r}{p}\\
\rightarrow i_{max,1}&=m+\lfloor \frac{r}{p}\rfloor .
\end{align}

If $n$ is the number of straight paths of the billiard ball trajectory we thus have for the circle division of the first $m$ reflections
\begin{equation}
f_m=\sum_{n=0}^{m+\lfloor \frac{r}{p}\rfloor} 1.
\end{equation}
Note that since $r<p$ the term $\lfloor \frac{r}{p}\rfloor$ may be omitted. Nevertheless $\lfloor \frac{r}{p}\rfloor$ will appear in the following terms again and will play a crucial role for the general formula, thus we keep it at this point in the sum. The sequence of circle division regions $f_n$ generated by the trajectory for the first $m$ scatterings thus is given by (see for a graphical representation figure \ref{fig:example1} \textit{left})
\begin{equation}
f_n: 1\underbrace{\xrightarrow{+1} 2 \xrightarrow{+1} 3 \xrightarrow{+1} \dots \xrightarrow{+1}  }_{m-times}m+1.
\end{equation}

\begin{figure}[htb]
\begin{minipage}{0.3\textwidth}
\centering
\begin{tikzpicture}[scale=1]
\draw (0,0) circle (2);
\foreach \i in {0,1,2,3} \draw[blue] (360*3/14*\i:2)--({360*3/14*(\i+1)}:2);

%\draw[blue] (0,0) circle (0.7137);
%\draw[red] (0,0) circle (0.4939);
\end{tikzpicture}
\end{minipage}
\hfill
\begin{minipage}{0.3\textwidth}
\centering
\begin{tikzpicture}[scale=1]
\draw (0,0) circle (2);
\foreach \i in {0,1,2,3} \draw[blue] (360*3/14*\i:2)--({360*3/14*(\i+1)}:2);
\foreach \i in {4} \draw[red] (360*3/14*\i:2)--({360*3/14*(\i+1)}:2);

%\draw[blue] (0,0) circle (0.7137);
%\draw[red] (0,0) circle (0.4939);
\end{tikzpicture}
\end{minipage}
\hfill
\begin{minipage}{0.3\textwidth}
\centering
\begin{tikzpicture}[scale=1]
\draw (0,0) circle (2);
\foreach \i in {0,1,2,3} \draw[blue] (360*3/14*\i:2)--({360*3/14*(\i+1)}:2);
\foreach \i in {4} \draw[red] (360*3/14*\i:2)--({360*3/14*(\i+1)}:2);
\foreach \i in {5,6,7,8} \draw[green!60!black] (360*3/14*\i:2)--({360*3/14*(\i+1)}:2);

\end{tikzpicture}
\end{minipage}
\caption{\textit{left:} First $m$ divisions for $\vartheta=\frac{3}{4\cdot 3+2}\cdot 2\pi$. ~~~\textit{middle:} first performance of a rotation greater than $2\pi$   \textit{right:} following $m-1+\lfloor \frac{2r}{p}\rfloor -\lfloor \frac{r}{p}\rfloor$ divisions for $2\pi < n\cdot \vartheta <4\pi $ \label{fig:example1}}
\end{figure}

In the next step the trajectory is going to cross the initial point $P_0$ for the first time and thus will cross the line generated by the first path (see figure \ref{fig:example1} \textit{middle}). As a direct consequence for the path with $n=m+1$ the circle is divided into two additional regions since one line is crossed. The number of divisions for the first $m+1$ steps thus reads
\begin{equation}
f_{m+1}=\left(\sum_{n=0}^{m+\lfloor \frac{r}{p}\rfloor} 1\right)+2,
\end{equation}
resulting in the sequence
\begin{equation}
f_n: 1\underbrace{\xrightarrow{+1} 2 \xrightarrow{+1} 3 \xrightarrow{+1} \dots \xrightarrow{+1} }_{m-times} m+1 \xrightarrow{+2} m+3.
\end{equation}

Having performed the first $m+1$ steps the following paths (see figure \ref{fig:example1} right) will always cross two lines since two consecutive scattering points will always inhabit one scattering point from the first $m$ paths. From crossing two lines one will generate three new regions each time. The number of reflections $i_{max,2}$ in this situation that can be performed in this way can be determined by the knowledge of the rotation angle performed in the first $(m+1)$ reflections, which then crossed the initial point $P_0$ once and the restriction, that the point $P_0$ isn't crossed twice, i.e. 
\begin{align}
i_{max,2}\cdot \vartheta +(m+1+\lfloor \frac{r}{p}\rfloor )\cdot \vartheta &<4\pi\\
\leftrightarrow i_{max,2}&< \frac{4\pi}{\vartheta}-m-1 -\lfloor \frac{r}{p}\rfloor = m-1-\lfloor \frac{r}{p}\rfloor+\frac{2r}{p}\\
\rightarrow i_{max,2}&=m-1 +\lfloor \frac{2r}{p}\rfloor -\lfloor \frac{r}{p}\rfloor.
\end{align}
For the total number of regions thus produced in the first $2m+\lfloor \frac{2r}{p}\rfloor$ thus follows
\begin{equation}
f_{2m+\lfloor \frac{2r}{p}\rfloor}=\left(\sum_{n=0}^{m+\lfloor \frac{r}{p}\rfloor} 1\right)+2+\left(\sum_{n=m+2+\lfloor \frac{r}{p}\rfloor}^{2m+\lfloor \frac{2r}{p}\rfloor} 3\right)
\end{equation}
and for the related region sequence
\begin{equation}
f_n: 1\underbrace{\xrightarrow{+1} 2 \xrightarrow{+1} 3 \xrightarrow{+1} \dots \xrightarrow{+1} }_{m-times} m+1 \xrightarrow{+2} m+3 \underbrace{\xrightarrow{+3} m+6 \xrightarrow{+3}  \dots \xrightarrow{+3} }_{m-1 +\lfloor \frac{2r}{p}\rfloor -\lfloor \frac{r}{p}\rfloor~times}m+3+3\cdot\left(m-1 +\lfloor \frac{2r}{p}\rfloor\right).
\end{equation}

We can extend this logic to a full periodic orbit. Each time when the trajectory crosses the initial point $P_0$ we obtain an even number of new regions since an odd number of lines is crossed. For the paths in between an even number of lines is crossed resulting in generating an odd number of new regions each. This pattern continues until the last round. The general formula for the number of regions thus is given by

$$ f_{total}=1+\left(\sum_{n=1}^{m+\lfloor \frac{r}{p}\rfloor} 1\right) + 2+\left(\sum_{n=m+2+\lfloor \frac{r}{p}\rfloor}^{2m+\lfloor \frac{2r}{p}\rfloor} 3\right) + 4+ \left(\sum_{n=2m+2+\lfloor \frac{2r}{p}\rfloor}^{3m+\lfloor \frac{3r}{p}\rfloor} 5\right) + \dots $$

\begin{equation}
 \dots + \left(\sum_{n=(p-1)m+2+\lfloor \frac{(p-1)r}{p}\rfloor}^{pm+r+1} 2p-1  \right)  \label{eq:general_formula}
\end{equation}

with the associated division sequence
\begin{equation*}
f_n: 1\underbrace{\xrightarrow{+1} \dots \xrightarrow{+1}  }_{m~times} m+1 \xrightarrow{+2} m+3 \underbrace{\xrightarrow{+3}  \dots \xrightarrow{+3} }_{m-1 +\lfloor \frac{2r}{p}\rfloor -\lfloor \frac{r}{p}\rfloor~times} m+3+3\cdot\left(m-1 +\lfloor \frac{2r}{p}\rfloor -\lfloor \frac{r}{p}\rfloor \right) \xrightarrow{+4} \dots
\end{equation*}

\begin{equation}
\dots \underbrace{\xrightarrow{+2p-1}  \dots \xrightarrow{+2p-1} }_{m-1 +\lfloor \frac{pr}{p}\rfloor -\lfloor \frac{(p-1)r}{p}\rfloor~times} p\cdot(mp+r)+1 =pq+1.
\end{equation}

Note that the last sum in \eqref{eq:general_formula} runs to $pm+r+1$ instead of $pm+r$, i.e. there is one additional summand $2p-1$ resulting from the fact that the very last path ends exactly on the initial point $P_0$ and thus it crosses the same number of lines as the paths before.

As a small remark concerning the floor function holds $\lfloor \frac{k\cdot r}{p}\rfloor \in\{0;\dots ;r-1\} ~\forall k\in\{1;\dots;p-1\}$.

An interesting simplification takes place when for the rest term holds $r=1$. In these cases the sequence of the number of regions is simplified by the fact that
\begin{equation}
\lfloor \frac{k}{p}\rfloor =0~\forall k\in\{1;\dots ;p-1\}.
\end{equation}
The simplified version of \eqref{eq:general_formula} thus reads
\begin{equation}
f_{total}(r=1)=1+\left(\sum_{n=1}^m 1\right)+2+\left(\sum_{n=m+2}^{2m}3\right)+\dots +\left(\sum_{n=(p-1)m+2}^{pm+2} 2p-1\right)
\end{equation}
resulting in an integer sequence of alternating $(m-1)$ summations except the first $m$ and last $m$ reflections

\begin{equation}
f_n: 1\underbrace{\xrightarrow{+1} \dots \xrightarrow{+1}  }_{m~times} m+1 \xrightarrow{+2} m+3 \underbrace{\xrightarrow{+3}  \dots \xrightarrow{+3} }_{(m-1)~times} 4m \xrightarrow{+4} 4m+4 \underbrace{\xrightarrow{+5}  \dots \xrightarrow{+5} }_{(m-1)~times} \dots \underbrace{\xrightarrow{+2p-1}  \dots \xrightarrow{+2p-1} }_{m~times} pq+1.
\end{equation}

\begin{figure}[h!]
\centering

\begin{tikzpicture}[scale=0.8]
\draw (0,0) circle (2);
\foreach \i in {0,1,2,3,4,5,6,7,8,9,10,11,12} \draw (360*3/13*\i:2)--({360*3/13*(\i+1)}:2);
\end{tikzpicture}

$f_{13}=40$

\caption{Graphical representation for the division sequence in the case $\vartheta=\frac{3}{13}\cdot 2\pi$. \label{fig:example3}}
\end{figure}
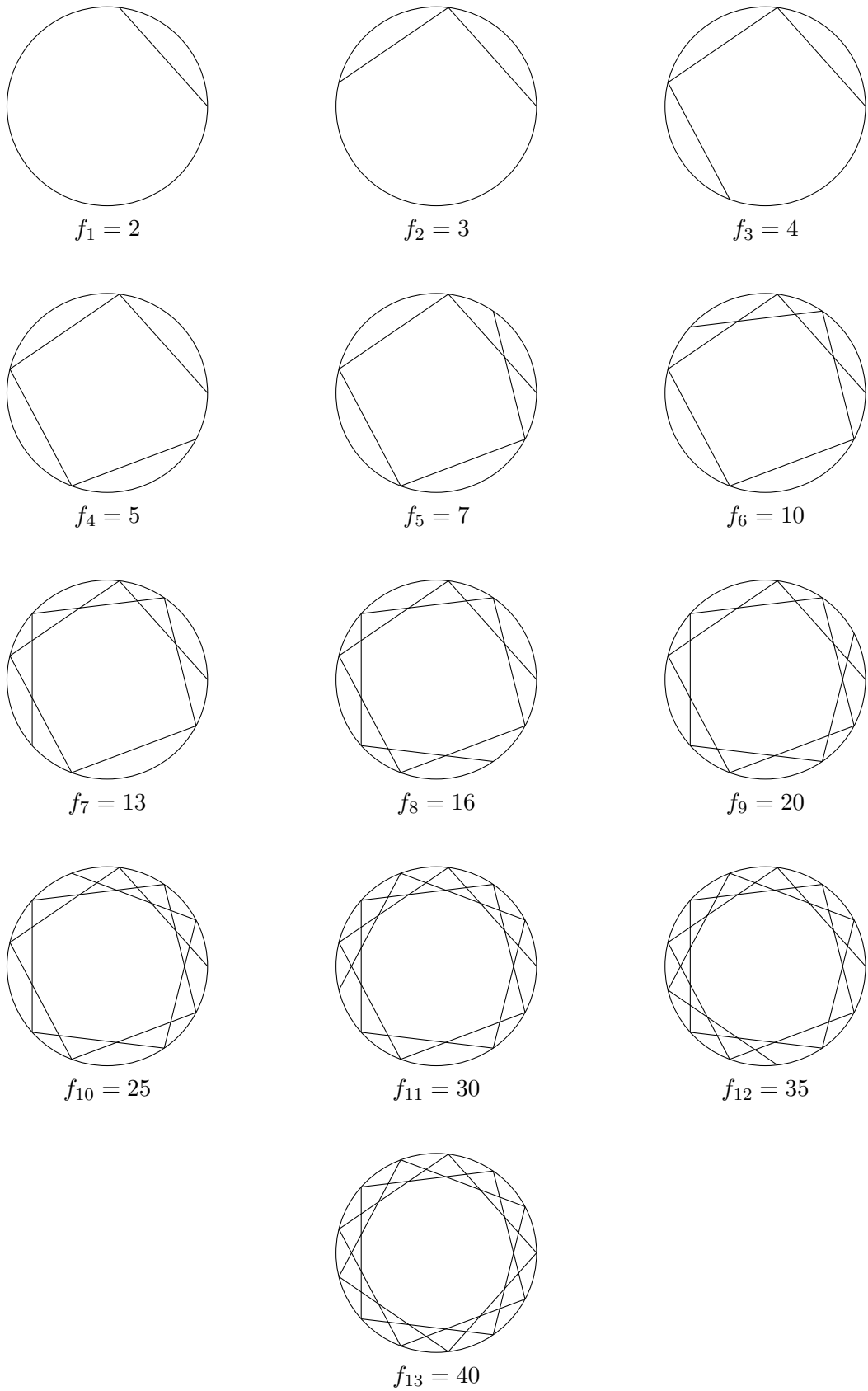

\vspace*{2cm}
%\newpage
To illustrate the last formula one should consider the example $\vartheta=\frac{3}{13}\cdot 2\pi$, i.e. $(p,m,r)=(3,4,1)$ shown in figure \ref{fig:example3}. The sequence for the regions reads
\begin{equation}
f_n:1\underbrace{\rightarrow 2\rightarrow 3\rightarrow 4 \rightarrow}_{4~ times~ +1} 5 \rightarrow 7 \underbrace{\rightarrow 10 \rightarrow 13 \rightarrow}_{3 ~times~+3}16 \rightarrow 20 \underbrace{\rightarrow 25 \rightarrow 30 \rightarrow 35 \rightarrow}_{4~times~+5} 40.
\end{equation}

\newpage
\begin{figure}[htb]
\begin{minipage}{0.3\textwidth}
\centering
\begin{tikzpicture}[scale=0.8]
\draw (0,0) circle (2);
\foreach \i in {0} \draw (360*3/13*\i:2)--({360*3/13*(\i+1)}:2);
\end{tikzpicture}

$f_1=2$
\end{minipage}
\hfill
\begin{minipage}{0.3\textwidth}
\centering
\begin{tikzpicture}[scale=0.8]
\draw (0,0) circle (2);
\foreach \i in {0,1} \draw (360*3/13*\i:2)--({360*3/13*(\i+1)}:2);
\end{tikzpicture}

$f_2=3$
\end{minipage}
\hfill
\begin{minipage}{0.3\textwidth}
\centering
\begin{tikzpicture}[scale=0.8]
\draw (0,0) circle (2);
\foreach \i in {0,1,2} \draw (360*3/13*\i:2)--({360*3/13*(\i+1)}:2);
\end{tikzpicture}

$f_3=4$
\end{minipage}
\end{figure}

\begin{figure}[htb]
\begin{minipage}{0.3\textwidth}
\centering
\begin{tikzpicture}[scale=0.8]
\draw (0,0) circle (2);
\foreach \i in {0,1,2,3} \draw (360*3/13*\i:2)--({360*3/13*(\i+1)}:2);
\end{tikzpicture}

$f_4=5$
\end{minipage}
\hfill
\begin{minipage}{0.3\textwidth}
\centering
\begin{tikzpicture}[scale=0.8]
\draw (0,0) circle (2);
\foreach \i in {0,1,2,3,4} \draw (360*3/13*\i:2)--({360*3/13*(\i+1)}:2);
\end{tikzpicture}

$f_5=7$
\end{minipage}
\hfill
\begin{minipage}{0.3\textwidth}
\centering
\begin{tikzpicture}[scale=0.8]
\draw (0,0) circle (2);
\foreach \i in {0,1,2,3,4,5} \draw (360*3/13*\i:2)--({360*3/13*(\i+1)}:2);
\end{tikzpicture}

$f_6=10$
\end{minipage}
\end{figure}

\begin{figure}[htb]
\begin{minipage}{0.3\textwidth}
\centering
\begin{tikzpicture}[scale=0.8]
\draw (0,0) circle (2);
\foreach \i in {0,1,2,3,4,5,6} \draw (360*3/13*\i:2)--({360*3/13*(\i+1)}:2);
\end{tikzpicture}

$f_7=13$
\end{minipage}
\hfill
\begin{minipage}{0.3\textwidth}
\centering
\begin{tikzpicture}[scale=0.8]
\draw (0,0) circle (2);
\foreach \i in {0,1,2,3,4,5,6,7} \draw (360*3/13*\i:2)--({360*3/13*(\i+1)}:2);
\end{tikzpicture}

$f_8=16$
\end{minipage}
\hfill
\begin{minipage}{0.3\textwidth}
\centering
\begin{tikzpicture}[scale=0.8]
\draw (0,0) circle (2);
\foreach \i in {0,1,2,3,4,5,6,7,8} \draw (360*3/13*\i:2)--({360*3/13*(\i+1)}:2);
\end{tikzpicture}

$f_9=20$
\end{minipage}
\end{figure}

\begin{figure}[h!]
\begin{minipage}{0.3\textwidth}
\centering
\begin{tikzpicture}[scale=0.8]
\draw (0,0) circle (2);
\foreach \i in {0,1,2,3,4,5,6,7,8,9} \draw (360*3/13*\i:2)--({360*3/13*(\i+1)}:2);
\end{tikzpicture}

$f_{10}=25$
\end{minipage}
\hfill
\begin{minipage}{0.3\textwidth}
\centering
\begin{tikzpicture}[scale=0.8]
\draw (0,0) circle (2);
\foreach \i in {0,1,2,3,4,5,6,7,8,9,10} \draw (360*3/13*\i:2)--({360*3/13*(\i+1)}:2);
\end{tikzpicture}

$f_{11}=30$
\end{minipage}
\hfill
\begin{minipage}{0.3\textwidth}
\centering
\begin{tikzpicture}[scale=0.8]
\draw (0,0) circle (2);
\foreach \i in {0,1,2,3,4,5,6,7,8,9,10,11} \draw (360*3/13*\i:2)--({360*3/13*(\i+1)}:2);
\end{tikzpicture}

$f_{12}=35$
\end{minipage}
\end{figure}

\newpage
\section{Discussion of the special case $\vartheta =\frac{p}{2p+1}\cdot 2\pi$}
Considering the case $(m,r)=(2,1)$ corresponding to angles $\vartheta=\frac{p}{2p+1}\cdot 2\pi$ formula \eqref{eq:general_formula} for the number of regions $f_n$ significantly simplifies to 
\begin{equation}
f_n =1+1+1+2+3+4+\dots  + 2p-2 +2p-1 +2p-1,
\end{equation}
If $n$ is the number of performed scatterings assuming we start with $0\leq n \leq 2p+1$ the above sequence can be rewritten in terms of a formula for $f_n$ depending on $n$ as

\begin{equation}
\boxed{f_n=2-\delta_{n,0}+\sum_{i=1}^{n-1} i -\delta_{2p+1,n}=2-1\cdot(\delta_{n,0}+\delta_{n,2p+1})+\frac{n\cdot (n-1)}{2}} \label{eq:main}
\end{equation}

This formula tells us that despite of the first scattering ($n=1$) and the last scattering ($n=2p+1$) the number of circle divisions performed by the billiard ball trajectory is related to integer sequence of arithmetic series 
\begin{equation}
f_{n+1}=f_n+(n-1)
\end{equation}
which starts with $f_0=1$ and $f=1=2$. Further $\delta_{i,j}$ represents the Kronecker-Delta defined by

\begin{equation}
\delta_{i,j}=\begin{cases} 1 &if~i=j\\
0 & otherwise
\end{cases}
\end{equation}

Formula \eqref{eq:main} states our main (applicable) result of this paper and extends the already known examples of circle divisions either by lines \cite{line1,line2} or by chords \cite{chords1,chords2,chords3}.  

In particular, if $n=2p+1$ scatterings are performed, we again obtain the total number of regions that the circle has been divided in by the periodic trajectory
\begin{equation}
f_{total}=f_{2p+1}=2-1+\frac{(2p+1)\cdot (2p+1-1)}{2}=pq+1,
\end{equation}
where again $2p+1$ was defined as $q$ and thus matching with equation \eqref{eq:areas}.

\begin{figure}[htb]
\begin{minipage}{0.18\textwidth}
\centering
\begin{tikzpicture}[scale=0.7]
\draw (0,0) circle (2);
\foreach \i in {0,1,2,3} \draw (360*3/7*\i:2)--({360*3/7*(\i+1)}:2);
\end{tikzpicture}
$f_4=8$
\end{minipage}
\hfill
\begin{minipage}{0.18\textwidth}
\centering
\begin{tikzpicture}[scale=0.7]
\draw (0,0) circle (2);
\foreach \i in {0,1,2,3,4} \draw (360*3/7*\i:2)--({360*3/7*(\i+1)}:2);
\end{tikzpicture}
$f_5=12$
\end{minipage}
\hfill
\begin{minipage}{0.18\textwidth}
\centering
\begin{tikzpicture}[scale=0.7]
\draw (0,0) circle (2);
\foreach \i in {0,1,2,3,4,5} \draw (360*3/7*\i:2)--({360*3/7*(\i+1)}:2);
\end{tikzpicture}
$f_6=17$
\end{minipage}
\hfill
\begin{minipage}{0.18\textwidth}
\centering
\begin{tikzpicture}[scale=0.7]
\draw (0,0) circle (2);
\foreach \i in {0,1,2,3,4,5,6} \draw (360*3/7*\i:2)--({360*3/7*(\i+1)}:2);
\end{tikzpicture}
$f_7=22$
\end{minipage}
\caption{Billiard paths and number of circle divisions for $\vartheta =\frac{3}{7}\cdot 2\pi$. \label{fig:example2}}
\end{figure}
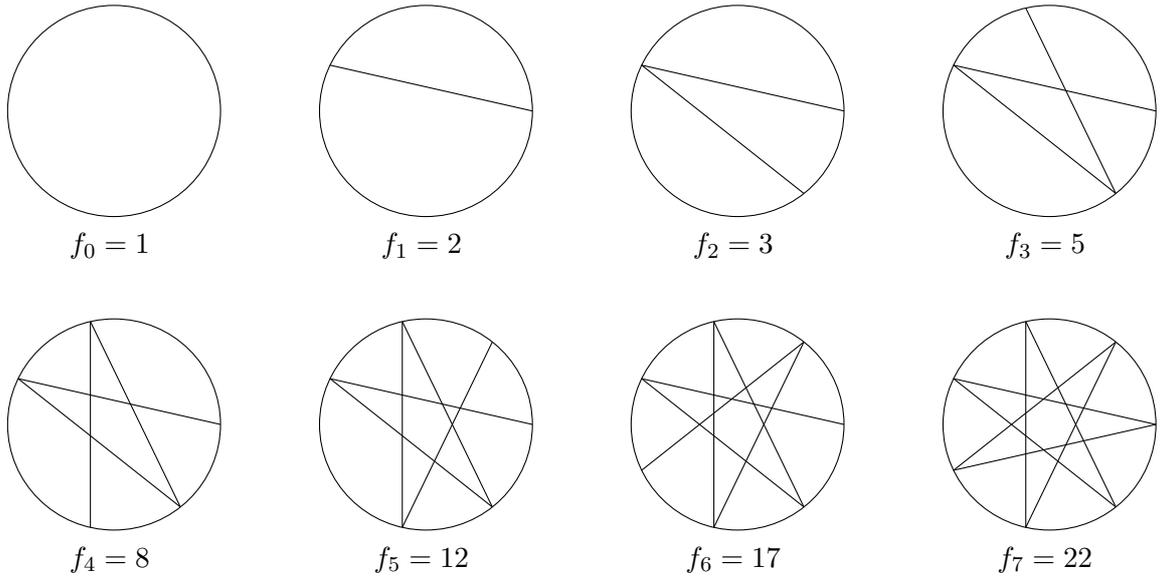

As an example the sequence for the number of regions for circle divisions of the example in figure \ref{fig:example2} is given by
\begin{equation}
f_n: 1\xrightarrow{+1} 2 \xrightarrow{+1} 3 \xrightarrow{+2} 5 \xrightarrow{+3} 8 \xrightarrow{+4} 12 \xrightarrow{+5} 17 \xrightarrow{+6-1} 22
\end{equation}

\begin{figure}[htb]
\begin{minipage}{0.18\textwidth}
\centering
\begin{tikzpicture}[scale=0.7]
\draw (0,0) circle (2);
%\foreach \i in {0,1,2,3} \draw (360*3/7*\i:2)--({360*3/7*(\i+1)}:2);
\end{tikzpicture}
$f_0=1$
\end{minipage}
\hfill
\begin{minipage}{0.18\textwidth}
\centering
\begin{tikzpicture}[scale=0.7]
\draw (0,0) circle (2);
\foreach \i in {0} \draw (360*3/7*\i:2)--({360*3/7*(\i+1)}:2);
\end{tikzpicture}
$f_1=2$
\end{minipage}
\hfill
\begin{minipage}{0.18\textwidth}
\centering
\begin{tikzpicture}[scale=0.7]
\draw (0,0) circle (2);
\foreach \i in {0,1} \draw (360*3/7*\i:2)--({360*3/7*(\i+1)}:2);
\end{tikzpicture}
$f_2=3$
\end{minipage}
\hfill
\begin{minipage}{0.18\textwidth}
\centering
\begin{tikzpicture}[scale=0.7]
\draw (0,0) circle (2);
\foreach \i in {0,1,2} \draw (360*3/7*\i:2)--({360*3/7*(\i+1)}:2);
\end{tikzpicture}
$f_3=5$
\end{minipage}
\end{figure}

It is remarkable that the integer sequence
\begin{equation}
\sum_i i,
\end{equation}
which is known in many different areas of pure mathematics and physics also appears in the circle divisions of billiard trajectories for each angle $\vartheta =\frac{p}{2p+1}\cdot 2\pi,~p\in\mathbb{N}$, now.

\section{Conclusion and Outlook}\label{sec:conclusion}
In this paper we have provided a general formula for the region division sequence of a rational billiard path inside a circle. By this we have analyzed another examples of circle divisions besides the already known cases of divisions  produced by chords or by lines. We have presented a simple formula for the number of regions following after each scattering $n$ for certain angles $\vartheta =\frac{p}{2p+1}\cdot 2\pi$, $p\in \mathbb{N}$. In particular, this formula containes a arithmetic series.

For future works it will be interesting to see if similar sequences arise from elliptical, triangular or general polygonal domains.

%%%%%%%%%%%%%%%%%%%%%%%%%%%%%%%%%%%%%%%%%%%%%%%%%%%%%%%%%%%%%%%%%%%%%%%%%%%%%%%%%%%%%%%%%%%%%%%%%%%%%%%%%%%%%%%%%%%%%%
%\section*{Acknowledgements}
%We would like to thank Jan Bickel for the useful discussions.
%%%%%%%%%%%%%%%%%%%%%%%%%%%%%%%%%%%%%%%%%%%%%%%%%%%%%%%%%%%%%%%%%%%%%%%%%%%%%%%%%%%%%%%%%%%%%%%%%%%%%%%%%%%%%%%%%%%%%%

%%%%%%%%%%%%%%%%%%%%%%%%%%%%%%%%%%%%%%%%%%%%%%%%%%%%%%%%%%%%%%%%%%%%%%%%%%%%%%%%%%%%%%%%%%%%%%%%%%%%%%%%%%%%%%%%%%%%%%
\appendix

%\section{Literature}

\vspace*{8cm}
\begin{center}
\textit{A special thanks to Jan Bickel for a useful discussion.}
\end{center}

\end{document}